\documentclass[11pt]{amsart}
\usepackage{amssymb,amsmath,amsthm}
\usepackage{epsfig}
\newtheorem{theorem}{Theorem}[section]
\newtheorem{lemma}[theorem]{Lemma}

\newtheorem{corollary}[theorem]{Corollary}

\theoremstyle{remark}
\newtheorem{definition}[theorem]{Definition}
\newtheorem{remark}[theorem]{Remark}

\newcommand\Z{{\mathbb Z}}
\begin{document}

\title{The Quandle of the Trefoil Knot as the Dehn Quandle of the Torus}
\date{April 10, 2008}
\author{Maciej Niebrzydowski}
\address[Maciej Niebrzydowski]{Department of Mathematics\\
	 University of Louisiana at Lafayette\\
	 Lafayette, LA 70504}
\email{mniebrz@gmail.com}
\author{J\'ozef H. Przytycki}
\address[J\'ozef H. Przytycki]{Department of Mathematics\\
         The George Washington University\\
         Washington, D.C. 20052\\
         Department of Mathematics\\
	 University of Maryland\\
	 College Park, MD 20742}
\email{przytyck@gwu.edu}


\thispagestyle{empty}

\begin{abstract}
We prove that the fundamental quandle of the trefoil knot is isomorphic to
the projective primitive subquandle of transvections of the symplectic space
$\Z \oplus \Z$.
The last quandle can be identified with the Dehn quandle of the torus and
the cord quandle on a 2-sphere with four punctures. We also show that the
fundamental quandle of the long trefoil knot is isomorphic to the cord quandle
on a 2-sphere with a hole and three punctures.
\end{abstract}

\maketitle

\section{Definitions and preliminary facts}  \label{1}

\begin{definition}
A {\it quandle}, $X$, is a set with a binary operation 
$(a, b) \mapsto a * b$
such that
\begin{enumerate}
\item For any $a \in X$, 
$a* a =a$.

\item For any $a,b \in X$, there is a unique $c \in X$ such that 
$a= c*b$. 

\item For any $a,b,c \in X$,
$ (a*b)*c=(a*c)*(b*c)$ (right distributivity). 
\end{enumerate}
\end{definition}

Note that the second condition can be replaced with the following requirement:
the operation $*b\colon Q\to Q$, defined by $*b(x)=x*b$, is a bijection. The inverse map to $*b$ will be denoted by
$\overline{*}b$. 

\begin{definition}
A {\it rack} is a set with a binary operation that satisfies 
(2) and (3).
\end{definition}
According to \cite{F-R}, the earliest discussion on racks is due to J. Conway and G. Wraith, who studied racks in the context of the conjugacy operation in a group. They regarded a rack as the wreckage of a group left behind after the group operation is discarded and only the notion of conjugacy remains. 
The notion of quandle was introduced independently by D. Joyce \cite{Joy} and S. Matveev \cite{Mat}.

The fundamental quandle of the oriented knot is a classifying invariant of classical unoriented knots (see \cite{Joy} for details). Its generators correspond to the arcs of the knot diagram, and relations correspond to crossings. They are of the form: $x_i\ast x_k=x_j$ or $x_i\,\bar{\ast}\,x_k=x_j$, depending on the type of the crossing, where $x_i$ and $x_j$ are generators assigned to the under-arcs, and $x_k$ is assigned to the over-arc of the crossing.
Just like in the case of the fundamental groups, it is not easy to decide whether two given knot quandles are isomorphic.

The following are some of the most common examples of quandles.
\begin{itemize}
\item[-]
Any group $G$ with conjugation 
as the quandle operation:\\ $a*b=b^{-1} a b$. The quandle $(G,*)$ is denoted by $Conj(G)$.
\item[-]
Let $n$ be a positive integer.
For elements  $i, j \in \{ 0, 1, \ldots , n-1 \}$, define
$i\ast j \equiv 2j-i \pmod{n}$.
Then $\ast$ defines a quandle
structure  called the {\it dihedral quandle},
  $R_n$.
It can be identified with  the
set of reflections of a regular $n$-gon
with conjugation
as the quandle operation.
\item[-]
Any $\mathbb{Z}[t, t^{-1}]$-module $M$
is a quandle with
$a*b=ta+(1-t)b$, for $a,b \in M$, called an {\it  Alexander  quandle}.
Moreover, if $n$ is a positive integer, then
$\mathbb{Z}_n[t, t^{-1}]/(h(t))$
is a quandle for
a Laurent polynomial $h(t)$.
\end{itemize}

The last example can be vastly generalized \cite{Joy}; for any group $G$ and its automorphism $\tau\colon G\to G$, $G$ becomes a quandle when equipped with the operation $g*h=\tau(gh^{-1})h$. If we consider the anti-automorphism $\tau(g)=g^{-1}$, we obtain another well known quandle, $Core(G)$, with
$g*h=hg^{-1}h$. 

\section{Dehn quandles and symplectic quandles}\label{Section 2} 

In this section we recall (after J.~Zablow \cite{Za-1,Za-2} and D.~Yetter \cite{Ye-1,Ye-2}; compare also \cite{K-M}) 
the concept of Dehn quandle of an orientable surface, and related definition of a symplectic quandle.

Let $F$ be an orientable surface and let $C(F)$ denote the isotopy classes of simple closed curves in $F$.
For any curve $c\in C(F)$, we consider the positive (right-handed) or negative (left-handed) Dehn twist about this curve, denoted as $t^+_c$ and $t^-_c$ respectively.

The following facts are needed when defining the Dehn quandle of a surface $F$.
\begin{enumerate}
\item[-] $t^+_c$ fixes the curve $c$ up to isotopy;
\item[-] Positive and negative Dehn twists about the same curve are inverse to each other up to isotopy;
\item[-] Positive Dehn twists along isotopic simple closed curves are isotopic as diffeomorphisms;
\item[-] The images of simple closed curves under isotopic Dehn twists are isotopic. 
\end{enumerate}

Thus, it makes sense to consider the following definition, in which the same symbol denotes the isotopy class and a representative curve.

\begin{definition}
The {\it Dehn quandle}, $Dehn(F)$, of an orientable surface $F$, is the set of isotopy classes of simple closed curves in $F$, equipped with the operations
$$x*y=t^+_y(x),$$
$$x\ \bar*\ y=t^-_y(x).$$
\end{definition}

For a detailed proof that $Dehn(F)$ satisfies quandle axioms see \cite{Za-1}.

The Dehn twist $t_y^+$ acts on $H_1(F)$ in a natural way. The action depends only on the homology class of $y$, and preserves the intersection form on $H_1(F)$. This motivates the following definition.

\begin{definition}
Let $R$ be a commutative ring, $M$ be a module over $R$, and let $<\ ,\ >\colon M\times M\to R$ be a bilinear form. Consider the operation $x*y=x-<x,y>y$.
Then: 
\begin{enumerate}
\item[(i)] If $<,>$ is antisymmetric, i.e., $<x,y>=-<y,x>$, then $(M,*)$ is a rack.
\item[(ii)] If $<,>$ satisfies $<x,x>=0$ for all $x$, then $(M,*)$ is a quandle\footnote{$<x,x>$ implies $<x,y>=-<y,x>$, and the inverse holds if $2$ is not a zero divisor in $R$.}.
\end{enumerate}
We refer to above quandle as {\it symplectic} quandle if M is free and the form $<\ ,\ >$ is non-degenerate.  
\end{definition}

The structure of symplectic quandles was recently studied in \cite{N-N}. 

\section{Proof of the main theorem}\label{Section 3}

\begin{figure}
\begin{center}
\includegraphics[height=4cm]{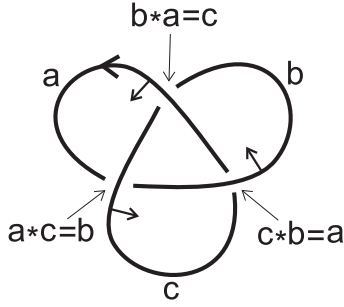}
\caption{Generators and relations of the fundamental quandle of the trefoil knot.\label{trefoil} }
\end{center}
\end{figure}

By definition, the fundamental quandle of the trefoil knot, $Q(3_1)$, has
presentation 
$$\{a,b,c\ |\ a*c=a, b*a=c, c*b=a \}=\{a,b\ |\ a*b*a=b, b*a*b=a \};$$ compare Fig.\ref{trefoil}. An important property of $Q(3_1)$ is that it satisfies the braid type relation:
\begin{equation}\label{braid}
x*a*b*a = x*b*a*b, 
\end{equation}
for any $x \in Q(3_1)$.
In particular, it allows 
a homomorphism from the 3-braid group, $B_3$, to the group of inner automorphisms of
$Q(3_1)$, sending $\sigma_1$ to $*a$ and 
$\sigma_2$ to $*b$, where $\sigma_1$, $\sigma_2$ are standard generators of $B_3$.

Let us first prove relation (\ref{braid}). 
The equation $a*b*a=b$ is equivalent to $a=b\bar*a\bar*b,$ and it follows that
$x*a=x*(b\bar*a\bar*b)$, for any $x \in Q(3_1)$. Thus, $x*a=x*b*a*b\bar*a\bar*b$, and after applying $*b*a$ to both sides of the last equation, we get the required relation
$x*a*b*a = x*b*a*b$. We remark that we used only relation $a*b*a=b$ to get (\ref{braid}). It follows that above is true also for the fundamental quandle of the long trefoil knot, in which the relation $b*a*b=a$ does not hold.

The following theorem suggests that there could be a strong connection between the fundamental knot quandles and Dehn quandles.

\begin{theorem}
The fundamental quandle of the trefoil knot, $Q(3_1)$, is isomorphic to the Dehn quandle of the torus, $Dehn(T^2)$. 
\end{theorem}

\begin{proof}

From the fact that $$\pi_1(T^2)=H_1(T^2)=\Z\oplus \Z,$$ follows that isotopy classes of unoriented curves 
on $T^2$ (and Dehn twists assigned to them) can be identified with relatively prime pairs $$\pm (a,b) \in \Z\oplus \Z/\pm 1,$$ in the space of orbits of the action of the multiplicative group $\{1, -1\}$ on 
$\Z\oplus \Z$ by scalar multiplication; in other words, with fractions 
$\frac{a}{b}\in Q\cup \frac{1}{0}$.
Furthermore, the action $*\beta$ 
given by the Dehn twist corresponding to the curve with ``slope" $\beta=\frac{c}{d}$ is a transvection, 
that is, for $\alpha=\frac{a}{b}$, we have $\alpha * \beta = \frac{a - Dc}{b - Dd},$ where $D$ is the determinant,
$D= ad-bc$. Indeed, one can easily check that $*\beta$ is 
an element of $PSL(2,\Z)$ given by the matrix: 
\begin{equation*} 
\left[
\begin{array}{cc}
1-dc & c^2 \\
-d^2 & 1+dc
\end{array}
\right],
\end{equation*}
and then 
\begin{equation*}
\left[
\begin{array}{cc}
1-dc & c^2 \\
-d^2 & 1+dc
\end{array}
\right] \left[
\begin{array}{c}
a\\
b\end{array}
\right] =
 \left[\begin{array}{c}
a-cD\\
b-dD\end{array} \right].
\end{equation*}
In particular, for $(c,d)=(1,0),$ we have 
\begin{equation*} *_{\beta} =
\left[
\begin{array}{cc}
1 & 1 \\
0 & 1%
\end{array}
\right],
\end{equation*}

\noindent and the matrix for $(c,d)=(0,1)$ is

\begin{equation*} *_{\beta} =
\left[
\begin{array}{cc}
1 & 0 \\
-1 & 1
\end{array}
\right],
\end{equation*}
that is, we get generators of $PSL(2,\Z)$.

We remark that the determinant $D((a,b),(c,d))=ad-bc$ is a symplectic form on $\Z\oplus \Z$, and the symplectic quandle operation is given by
$$(a,b)*(c,d)=(a-cD, b-dD).$$ Since $Dehn(T^2)$ corresponds to relatively prime (primitive) pairs from\\ $\Z\oplus \Z/\pm 1,$ we can say that
it can be identified with the projective primitive subquandle of the symplectic quandle on $\Z\oplus \Z$. 

The second part of the proof (the correspondence between fractions and elements of $Q(3_1)$) follows from Theorem \ref{main}.
\end{proof}

\begin{theorem}\label{main}
$Q(3_1)$ is a quandle isomorphic to the quandle of fractions 
$(Q\cup \frac{1}{0},*)$, 
where $(a/b) * (c/d) = \frac{a - Dc}{b - Dd},$ where $D=ad-bc$ is the determinant.
The isomorphism  $\phi: Q(3_1)\to (Q\cup \frac{1}{0},*)$ is given by
$\phi (a) = \frac{0}{1}$, $\phi (b) = \frac{1}{0}$.
\end{theorem}

\begin{proof}
The map $\phi$ is a quandle homomorphism because
$$\phi(a) * \phi(b) * \phi(a) = \frac{0}{1} * \frac{1}{0} * \frac{0}{1} = \frac{1}{1} * \frac{0}{1} 
= \frac{1}{0} = \phi(b) = \phi(a*b*a),$$
and similarly, 
$$\phi(b) * \phi(a) * \phi(b) = \frac{1}{0} * \frac{0}{1} * \frac{1}{0} = \frac{1}{-1} * \frac{1}{0} = 
\frac{0}{-1} = \frac{0}{1} = \phi(a) = \phi(b*a*b).$$
 
In order to prove that $\phi$ is an epimorphism, we are going to represent rational numbers as continued fractions. 
We write $[k_1; k_2, k_3,\ldots, k_n]$ to denote the continued fraction
$$ k_1 + \frac{1}{k_2+\frac{1}{k_3+\ldots {+\frac{1}{k_n}}}}.$$
Let us recall the algorithm for expanding any rational number $r$ into a continued fraction.
Let $k_1=[r]$ be the greatest integer not exceeding $r$. It follows that $\delta=r - k_1<1$ and $\delta\geq 0$.
If $\delta = 0$, the algorithm ends. Otherwise, let $r_2=\frac{1}{\delta}$, $k_2=[r_2]$, and $\delta=r_2-k_2<1$. It is not difficult to show that after a finite number of such steps we obtain $\delta=0$ and the algorithm ends.

In the continued fraction representation of $\alpha$ that we obtain from the above algorithm, $k_1$ is an integer, $k_2,\ldots, k_n$ are positive integers, and $k_n>1$. In fact, we have the following lemma \cite{Khi}.

\begin{lemma}\label{rep}
The correspondence between
\begin{enumerate}
\item finite continued fractions $[k_1 ; k_2, k_3,\ldots, k_n]$ with an integer $k_1$, positive integers $k_2,\ldots,k_n$ and $k_n>1$
\item[] and
\item rational numbers
\end{enumerate}
is one-to-one.
\end{lemma}

We will use the convention that $u*w^k = u*w*w*...*w$ ($k$-times $*w$) for $k>0$
and $u*w^k = u\bar*w\bar*w\bar*...\bar*w$ ($-k$-times $\bar*w$) for $k<0$.

\begin{lemma}\label{epi}
Let $[k_1 ; k_2, k_3,\ldots, k_n]$ be a continued fraction satisfying the conditions of Lemma \ref{rep} corresponding to a rational number $\frac{p}{q}$.
If $n$ is odd, then $$\phi (a*b^{k_n}*a^{-k_{n-1}}*...*b^{k_1})=\frac{p}{q}.$$
If $n$ is even, then $$\phi (b*a^{-k_n}*b^{k_{n-1}}*...*b^{k_1})=\frac{p}{q}.$$
\end{lemma} 
\begin{proof}
In the proof, the following formulas will be useful.
$$\frac{p}{q}*\Big(\frac{s}{t}\Big)^k=\frac{p-kDs}{q-kDt},$$
where $k>0$ and $D=pt-sq$; 
$$\frac{p}{q}*\Big(\frac{s}{t}\Big)^k=\frac{p+kDs}{q+kDt},$$
where $k<0$ and $D$ is as above.
Below we perform the inductive step in the proof of the first formula (the proof of the second formula is analogous):
$$\frac{p}{q}*\Big(\frac{s}{t}\Big)^k=\frac{p}{q}*\Big(\frac{s}{t}\Big)^{k-1}*\frac{s}{t}=\frac{p-(k-1)Ds}{q-(k-1)Dt}*\frac{s}{t}=
\frac{p-(k-1)Ds-Ds}{q-(k-1)Dt-Dt},$$ because $(p-(k-1)Ds)t-s(q-(k-1)Dt)=D$.

It follows that for positive $k$:
\begin{enumerate}
\item[(i)]  $\frac{p}{q}*\Big(\frac{0}{1}\Big)^k=\frac{p}{q-kp}=\frac{1}{\frac{q-kp}{p}}=\frac{1}{-k +\frac{q}{p}}$
\item[(ii)]  $\frac{p}{q}\ \bar*\ \Big(\frac{0}{1}\Big)^k=\frac{p}{q+kp}=\frac{1}{\frac{q+kp}{p}}=\frac{1}{k+ \frac{q}{p}}$
\item[(iii)] $\frac{p}{q}*\Big(\frac{1}{0}\Big)^k=\frac{p+kq}{q}=k+ \frac{p}{q}$
\item[(iv)] $\frac{p}{q}\ \bar*\ \Big(\frac{1}{0}\Big)^k=\frac{p-kq}{q}=-k +\frac{p}{q}$
\end{enumerate}

We will use induction on $n$ to prove formulas from Lemma \ref{epi}.
Assume that $\frac{p}{q}=[k_1; k_2, k_3,\ldots, k_n]$, where $n$ is odd.
Then, by inductive assumption, $$\phi (a*b^{k_n}*a^{-k_{n-1}}*...*b^{k_3})=[k_3; k_4, k_5,\ldots, k_n]=r.$$ 
Thus, $$\phi (a*b^{k_n}*a^{-k_{n-1}}*...*b^{k_3}*a^{-k_2})=r*{\phi(a)}^{-k_2}=r*{\Big(\frac{0}{1}\Big)}^{-k_2}=\frac{1}{k_2+ \frac{1}{r}},$$
from (ii), and
$$\phi (a*b^{k_n}*a^{-k_{n-1}}*...*b^{k_3}*a^{-k_2}*b^{k_1})=\Big(\frac{1}{k_2+ \frac{1}{r}}\Big)*{\phi(b)}^{k_1}=$$$$
\Big(\frac{1}{k_2+ \frac{1}{r}}\Big)*{\Big(\frac{1}{0}\Big)}^{k_1}=k_1+\frac{1}{k_2+ \frac{1}{r}}=k_1 + \frac{1}{k_2+\frac{1}{k_3+\ldots {+\frac{1}{k_n}}}}.$$
In this case we use (iii) if $k_1>0$ and (iv) if $k_1<0$.
The proof of the second formula is similar. It follows that $\phi$ is an epimorphism.
\end{proof} 
The fact that $\phi$ is a monomorphism follows from the following lemma.

\begin{lemma}\label{mono}
Each element of $Q(3_1)$ can be uniquely written in one of the following forms.
\begin{enumerate}
\item $a$ or $b$ or $a*b$ or $b*a$;
\item $a*b^{k_n}*a^{-k_{n-1}}*...*b^{k_1}$, where $n$ is odd, $k_2,\ldots,k_n$ are positive integers and $k_n>1$;
\item $b*a^{-k_n}*b^{k_{n-1}}*...*b^{k_1}$, where $n$ is even, $k_2,\ldots,k_n$ are positive integers and $k_n>1$.
\end{enumerate}
\end{lemma}
\begin{proof}
The following (operator level) relations are consequences of the relations $a*b*a=b$ and $b*a*b=a$ that are satisfied in $Q(3_1)$.
To shorten the expressions in this proof, we will denote $*a$ as $a$, $\bar*a$ as $\bar a$, $*b$ as $b$, and $\bar*b$ as $\bar b$.
\begin{enumerate}
\item[(R1)] $\ldots bab \ldots = \ldots aba \ldots$
\item[(R2)] $\ldots \bar b\bar a\bar b \ldots = \ldots \bar a\bar b\bar a \ldots$
\item[(R3)] $\ldots ba\bar b \ldots = \ldots \bar aba \ldots$
\item[(R4)] $\ldots b\bar a\bar b \ldots = \ldots \bar a\bar ba \ldots$
\item[(R5)] $\ldots \bar bab \ldots = \ldots ab\bar a \ldots$
\item[(R6)] $\ldots \bar b\bar ab \ldots = \ldots a\bar b\bar a \ldots$
\end{enumerate}
If we assume that the word satisfies all conditions of Lemma \ref{mono}, except the condition $k_n>1$ (i.e., we assume $k_n=1$), then we can fix this situation using
the equalities below and, if necessary, the induction on the length of the word (i.e., the total number of a's and b's appearing in the word representing the given element of $Q(3_1)$).
\begin{enumerate}
\item $ab\bar a=b\bar a\bar a$
\item $b\bar ab=abb$
\item $b\bar a\bar b=ab\bar b=a$
\end{enumerate}

Let us also note that our condition for the parity of $n$ follows from the first quandle axiom, $xx=x$ (or $x\bar x=x$), and the assumption that the words in the second and third class end in $b^{k_1}$.

Now we use induction on the length of the word to prove that the remaining conditions are also achievable.
We assume that our claim holds for any word $w$ of length $m$. We can also assume that $w$ is in the second or third class (the first class is easy to deal with). Let us extend $w$ by $a$ or $\bar a$ or $b$ or $\bar b$. Adding $b$ or $\bar b$ at the end of $w$ will not spoil our normal form.
Also, adding $\bar a$ if $w$ ends in $b^{s}$, where $s$ is positive, is permittable, as we allow $k_1$ to be $0$. Let us consider the three remaining cases.
\begin{enumerate}
\item[Case 1.] Assume that $w$ ends in ${\bar b}^s$, where $s$ is positive, and we consider the word $wa$. 
We have the following (operator level) relation in which we use brackets $[\ ]$ to indicate on which part of the word our relations are used.
\begin{align*}
\ldots {\bar b}^s a&=\ldots{\bar b}^{s-1}a[\bar a\bar ba]=\ldots{\bar b}^{s-1}ab\bar a\bar b=\ldots{\bar b}^{s-2}a[\bar a\bar ba]b\bar a\bar b\\
&=\ldots{\bar b}^{s-2}ab\bar a\bar bb\bar a\bar b=\ldots{\bar b}^{s-2}ab\bar a\bar a\bar b=\ldots=\ldots{\bar b}^{s-i}ab{\bar a}^i\bar b=\ldots\\
&=\ldots\bar bab{\bar a}^{s-1}\bar b=\ldots a[\bar a\bar ba]b{\bar a}^{s-1}\bar b= \ldots ab\bar a\bar bb{\bar a}^{s-1}\bar b=\ldots ab{\bar a}^s\bar b
\end{align*}
The last $s+2$ letters of $wa$ are now in the normal form. The $a$ that precedes these letters will cancel with $\bar a$ that appears before ${\bar b}^s$ (if there are no such $\bar a$ it means that $w$ was of the form $a^{{\bar b}^s}$ and we can use the first axiom of quandle). 
\item[Case 2.] Assume that $w$ ends in $b^s$, where $s$ is positive, and we consider the word $wa$.
We are going to use the following relation.
\begin{align*}
\ldots \bar ab^s&=\ldots b[\bar b \bar ab]b^{s-1}=\ldots ba\bar b\bar ab^{s-1}=\ldots ba[\bar b\bar ab]b^{s-2}=
\ldots baa\bar b\bar ab^{s-2}=\ldots\\
&=\ldots b a^i\bar b\bar a b^{s-i}=\ldots=\ldots ba^{s-1}[\bar b\bar ab]=
\ldots ba^{s-1}a\bar b\bar a=\ldots ba^{s}\bar b\bar a
\end{align*}
We have:
$$wa=\ldots \bar ab^{s}a=\ldots ba^{s}\bar b\bar aa=\ldots ba^{s}\bar b.$$
The length of the last expression is $m+1$. By inductive assumption, the first $m$ letters can be transformed into the word in the normal form, and since the $m+1$-th letter is $b$, we are done.
\item[Case 3.] Assume that $w$ ends in ${\bar b}^s$, where $s$ is positive, and we consider the word $w\bar a$.
We will need the following relation.
\begin{align*}
\ldots a{\bar b}^s&=\ldots \bar b[ba\bar b]{\bar b}^{s-1}=\ldots\bar b\bar a[ba\bar b]{\bar b}^{s-2}=
\ldots \bar b\bar a\bar aba{\bar b}^{s-2}=\ldots=\bar b{\bar a}^{i}ba{\bar b}^{s-i}\\&=\ldots=
\ldots \bar b{\bar a}^{s-1}[ba\bar b]=\ldots \bar b{\bar a}^{s-1}\bar aba=\ldots \bar b{\bar a}^{s}ba
\end{align*}
We will also use the relation $$\ldots b{\bar a}^t=\ldots \bar a{\bar b}^tab$$ that follows from the first relation and the fact that relations (R1)-(R6) are symmetric with respect to $a$ and $b$.
$$w\bar a=\ldots b{\bar a}^{k_2}{\bar b}^s\bar a=\ldots b{\bar a}^{k_2}\bar a[a{\bar b}^s]\bar a=
\ldots b{\bar a}^{k_2+1}\bar b{\bar a}^{s}ba\bar a=\ldots b{\bar a}^{k_2+1}\bar b{\bar a}^{s}b$$
The last expression is still too long to use induction (it has $m+3$ letters). That is why we use the second relation.
$$ \ldots [b{\bar a}^{k_2+1}]\bar b{\bar a}^{s}b=\ldots \bar a{\bar b}^{k_2+1}ab\bar b {\bar a}^{s}b=\ldots \bar a{\bar b}^{k_2+1}{\bar a}^{s-1}b$$
Now the last word has $m+1$ letters and the last letter is $b$, so we can use induction to end the proof.
\end{enumerate}
The uniqueness follows from the Lemma \ref{rep} and the fact that the map $\phi$ is well defined.
\end{proof}
The proof of Lemma \ref{mono} ends the proof of Theorem \ref{main}.
\end{proof}

\begin{remark} In the proof of Theorem \ref{main}, we could have used Ryder's theorem stating that for a prime knot, its fundamental quandle can be embedded into conjugation quandle of the fundamental group \cite{Ry}, together with the fact that $\pi_1(3_1)=B_3$ is a central $\Z$-extension of $PSL(2,\Z)$. Our goal, however, was to prove the theorem using elementary steps, on the level of quandles, so that the correspondence between elements of the two quandles becomes explicit.
\end{remark}

\section{Structure of the fundamental quandle of the long trefoil knot}

Our goal, in this section, is to prove that the quandle of the long trefoil knot can be viewed as a quandle of cords on the sphere with a hole and three punctures. In order to prove our theorem, we use Eisermann's description of quandles of long knots \cite{Eis}. For the convenience of the reader, we recall the facts from \cite{Eis} that we are going to use.

\subsection{Eisermann's description of the fundamental quandles of long knots}

Let $\pi_1$ be the fundamental group of a closed prime knot $K$, $m$ its meridian, $\pi_1'$ its commutator subgroup, and let $Q$ denote the conjugacy class of $m$ in $\pi_1$ with conjugation as a quandle operation. From the work of Ryder \cite{Ry}, it is known that in the case of prime knots, such a $Q$ is isomorphic to the fundamental quandle of the knot $K$.  In \cite{Eis}, the author considers the set $$\tilde{Q}(\pi_1,m)=\{(x,g')\in \pi_1\times\pi_1'\ |\ x=m^{g'}\}.$$ Here, we use the exponential notation to denote the conjugation, $x^y=y^{-1}xy$; later such notation will also be used for an action of the fundamental group on a cord quandle.

$\tilde{Q}(\pi_1,m)$ becomes a connected quandle when equipped with the operations
$$(x,g')*(y,h')=(x^y, g'x^{-1}y)=(x^y, m^{-1}g'(h')^{-1}mh'),$$
and
$$(x,g')\ \bar{*}\ (y,h')=(x^y, g'xy^{-1})=(x^y, mg'(h')^{-1}m^{-1}h').$$
Such operations already appear in the work of Joyce \cite{Joy}.
The quandle $\tilde{Q}(\pi_1,m)$ turns out to be isomorphic to the fundamental quandle, $Q_L$, of the long knot obtained from $K$ by breaking it at some point and extending the endpoints to infinity \cite{Eis}.
The map $p\colon \tilde{Q}(\pi_1,m) \to Q$ given by $p(m^{g'},g')=m^{g'}$ is a covering in the following sense.
\begin{definition}
A surjective quandle homomorphism $p\colon \tilde{Q}\to Q$ is called a {\it covering} if $p(\tilde{x})=p(\tilde{y})$ implies 
$\tilde{a}*\tilde{x}=\tilde{a}*\tilde{y}$ for all $\tilde{a}$, $\tilde{x}$, $\tilde{y}\in\tilde{Q}$. 
\end{definition}
\noindent Using the terminology from \cite{Joy}, we can say that 
$\tilde{x}$ and $\tilde{y}$ are behaviorally equivalent, that is, they act in the same way as operators.

As stated in \cite{Eis}, covering transformations for $p\colon \tilde{Q}(\pi_1,m) \to Q$ are given by the left action of
$\Lambda=C(m)\cap \pi_1'=<\lambda>$, where $C(m)$ denotes the centralizer of $m$, and $\lambda$ is the longitude of $K$ (see \cite{B-Z} for details on computing $\lambda$). The action is defined by $\lambda\cdot(m^{g'},g')=(m^{g'},\lambda g')$, and $<\lambda>$ acts freely and transitively on each fibre $p^{-1}(m^{g'})$.
\begin{definition}
A {\it representation} of a quandle $Q$ on a group $G$ is a map $\phi\colon Q\to G$ such that $\phi(a*b)=\phi(b)^{-1}\phi(a)\phi(b)$ for all $a$, $b\in Q$.
The map $\rho\colon Q\to Inn(Q)$, sending each quandle element $q$ to the corresponding operator $*q$, is called the natural representation of $Q$.
An {\it augmentation} consists of a representation $\phi\colon Q\to G$ together with a group homomorphism $\alpha\colon G\to Inn(Q)$ such that $\alpha\phi=\rho$. If $G$ is generated by the image $\phi(Q)$ (as in the case of knot groups), then the action of $G$ on $Q$ given by $\alpha$ is uniquely determined by the representation $\phi$, so we can say for simplicity that $\phi\colon Q\to G$ is an augmentation.
\end{definition}
\noindent As proven in \cite{Eis}, $p\colon \tilde{Q}(\pi_1,m) \to Q\subset \pi_1$ gives an augmentation. Here, the fundamental group acts on 
$\tilde{Q}(\pi_1,m)$ by $$(m^{g'},g')^h:=(m^{g'h}, m^{-\epsilon(h)}g'h),$$ where $h\in\pi_1$ and $\epsilon\colon \pi_1\to\Z$ is a homomorphism sending each element $q\in Q$ to 1. Furthermore, this action is by inner automorphisms. 

\subsection{Quandle of the long trefoil knot as a cord quandle}

\begin{figure}
\begin{center}
\includegraphics[height=4cm]{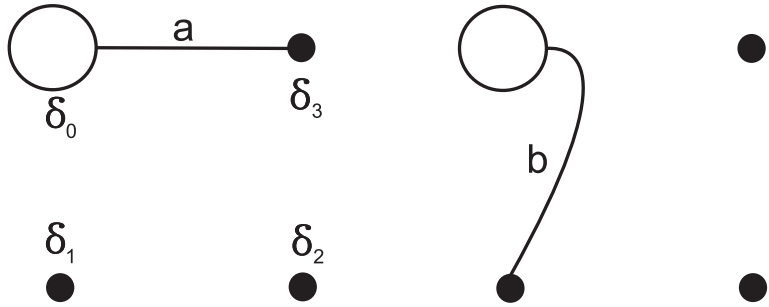}
\caption{Generators $a$ and $b$ of the cord quandle $\widehat{Q}$ .\label{gens} }
\end{center}
\end{figure}

We will define a quandle generated by the two cords $a$ and $b$ on a 2-sphere with four holes, $F_{0,4}$,
$\partial F_{0,4} = \{\delta_0, \delta_1, \delta_2, \delta_3\},$
see Fig.\ref{gens}.

\begin{figure}
\begin{center}
\includegraphics[height=4.6cm]{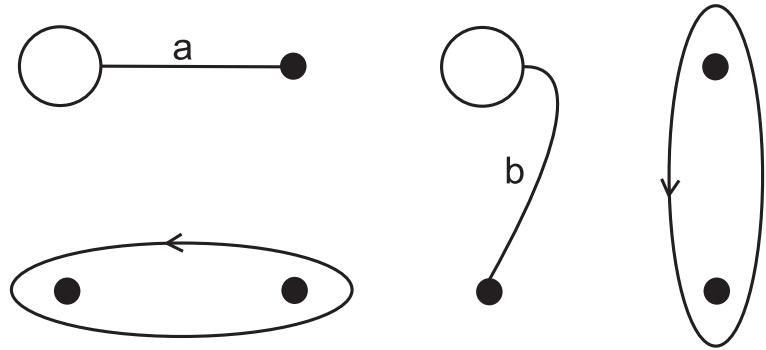}
\caption{Curves of Dehn half-twists $*a=\sigma_1$ and $*b=\sigma_2$.\label{dehn} }
\end{center}
\end{figure}

Consider the mapping class group of $F_{0,4}$ modulo the component
$\delta_0$; other components can be rotated and exchanged. With
this assumption, the mapping class group is the three-braid group,
$B_3= \{\sigma_1,\sigma_2\ |\ \sigma_1\sigma_2\sigma_1 =
\sigma_1\sigma_2\sigma_1 \},$ where $\sigma_1$ is the Dehn half-twist
(in the counter-clockwise direction) exchanging  $\delta_1$ and $\delta_2$ and
keeping $\delta_0$ fixed. Similarly, $\sigma_2$ is the Dehn half-twist
(in the counter-clockwise direction) exchanging  $\delta_2$ and $\delta_3$ and
keeping $\delta_0$ fixed; see Fig.\ref{dehn} and \cite{Bir}.

Now, consider the quandle $\widehat{Q}$ generated by the two arcs, $a$ and $b$,
with one endpoint at a fixed base point at $\delta_0$, as illustrated in the Figure \ref{gens}. 
That is, $\widehat{Q}$ consists of all arcs from the base point to
$\{\delta_1, \delta_2, \delta_3\}$,
with the convention that $*a = \sigma_1$ and $*b = \sigma_2$. The reader may wish to compare our definition with the definition of cord quandles given in \cite{K-M}.
The group of inner automorphisms of $\widehat{Q}$ is $B_3$. Since $B_3=\pi_1(3_1)$, the fundamental group of the trefoil knot acts on $\widehat{Q}$ by inner authomorphisms, and this action will be denoted using exponential notation.

\begin{figure}
\begin{center}
\includegraphics[height=4.6cm]{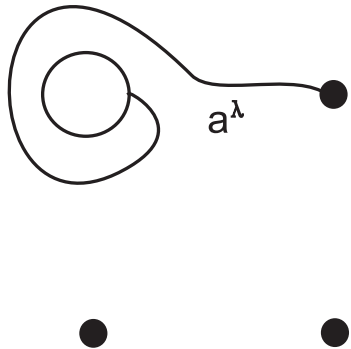}
\caption{$a^{\lambda}=a*a^{-4}*b*a*a*b$ -- full twist
along $\delta_0$ of the cord $a$.\label{fulltwist} }
\end{center}
\end{figure}

The most important identity in $\widehat{Q}$ is $a*b*a = b$. The element $a*a^{-4}*b*a*a*b$ is
obtained from $a$ by one clockwise twist along
$\delta_0$; see Fig.\ref{fulltwist}. Note that $\lambda=a^{-4}baab$ is a longitude for the trefoil knot.

We remark that given two elements $\alpha$, $\beta\in\widehat{Q}$, the quandle operation $\alpha*\beta$ can be realized as
$\alpha\tau[\beta]$, where $\tau[\beta]$ is the Dehn half-twist along the boundary of the regular neighborhood of $\beta$, exchanging the holes that are outside this neighborhood in the clockwise direction.

Our goal is the following theorem.
\begin{theorem} 
The quandle described above is isomorphic to the fundamental quandle of the long trefoil knot.
\end{theorem}

\begin{proof}
In our proof, $\pi_1$ denotes the fundamental group of the trefoil knot, and $\pi_1'$ denotes its commutator subgroup.
Recall that $\pi_1=B_3$ acts on $\tilde{Q}(\pi_1,m)$ by
\[(x,g')^g:=(x^g,g'x^{-\epsilon(g)}g)=(x^g, m^{-\epsilon(g)}g'g),\]
where $x=m^{g'}$, and it acts on $\widehat{Q}$ as a group of inner automorphisms.
Both quandles $\tilde{Q}(\pi_1,m)$ and $\widehat{Q}$ are connected. Connectedness of $\widehat{Q}$ can be seen from
the relation $a*b*a=b$, and the fact that $a$ and $b$ are the only generators. Connectedness of $\tilde{Q}(\pi_1,m)$
follows from the transitivity of the action of $\pi_1'$ and the fact that such action can be viewed as an action
by inner automorphisms \cite{Eis}.

Define a map $\Psi\colon \tilde{Q}(\pi_1,m)\to \widehat{Q}$ by
\[\Psi(m^{g'},g')=a^{g'}.\]
In particular, $\Psi(m,1)=a$, and $(m^{g'}, \lambda^k g')$ is sent to a cord wrapped around the hole $\delta_0$ $k$ times (in the clockwise direction if $k$ is positive, and counterclockwise otherwise), followed by a piece corresponding to $a^{g'}$.

We can also define a map in the opposite direction, $\Phi\colon\widehat{Q}\to \tilde{Q}(\pi_1,m)$, as follows.
Every element $q\in \widehat{Q}$ can be written as $a^w$, where $w\in \pi_1$.  
Because $\pi_1/\pi_1'$ is generated by the image of the meridian $m$, $\pi_1=<m>\pi_1'$, and every $w\in \pi_1$ can be uniquely written in the form $w=m^jw'$, where $j\in\Z$ and $w'\in \pi_1'$. 
Now, assume that $a^{w_1}=a^{w_2}$, where $w_1=m^iw_1'$, $w_2=m^jw_2'$, and $w_1'\neq w_2'$. The image of the element $a^{w_1}=a^{w_2}$ in $\pi_1$
equals $(w_1')^{-1}mw_1'=(w_2')^{-1}mw_2'$, so elements $(m^{w_1'},w_1')$, $(m^{w_2'},w_2')\in \tilde{Q}(\pi_1,m)$ are in the same fiber of the covering
$p\colon \tilde{Q}(\pi_1,m)\to Q(3_1)$. As noted in \cite{Eis}, the group $<\lambda>$ acts transitively on each fiber (the action is given by
$\lambda(m^{w'},w')=(m^{w'},\lambda w'))$. Therefore, $w_2'=\lambda^k w_1'$, for some nonzero integer $k$. It follows that $a^{w_1}=a^{m^iw_1'}$ cannot represent the same cord as $a^{w_2}=a^{m^jw_2'}$ because they differ by exactly $k$ twists around $\delta_0$.
Therefore, the map $\Phi\colon\widehat{Q}\to \tilde{Q}(\pi_1,m)$ given by
$$\Phi(a^w)=(m^{w'},w'),$$
where $w'\in \pi_1'$ is as above, is well defined.

We check that it is a homomorphism:
\begin{align*}
\Phi(a^{w_1}*a^{w_2})&=\Phi(a^{m^i w_1'm^{m^j w_2'}})=\Phi(a^{m^i w_1'(w_2')^{-1}m w_2'})=\Phi(a^{m^{i+1}m^{-1} w_1'(w_2')^{-1}m w_2'})\\
&=(m^{w_1'(w_2')^{-1}m w_2'}, m^{-1} w_1'(w_2')^{-1}m w_2'),
\end{align*}
where $w_1=m^iw_1'$, $w_2=m^jw_2'$, for some $w_1'$, $w_2'\in \pi_1'$, and $m^{-1} w_1'(w_2')^{-1}m w_2'\in\pi_1'$ because the sum of exponents in this word is zero. 
\[
\Phi(a^{w_1})*\Phi(a^{w_2})=(m^{w_1'},w_1')*(m^{w_2'},w_2')=(m^{w_1'm^{w_2'}}, m^{-1} w_1'(w_2')^{-1}m w_2'),
\]
as required.

We will show that both maps $\Psi$ and $\Phi$ are equivariant with respect to the action of $\pi_1$.
Let $x\in \pi_1$.
\[\Psi((m^{g'},g')^x)=\Psi(m^{g'x}, m^{-\epsilon(x)}g'x)=a^{m^{-\epsilon(x)}g'x}=a^{g'x}=(\Psi(m^{g'},g'))^x.\]
For the equivariance of the map $\Phi$, it is enough to consider the case when $x$ is an image of the quandle element $q=a^w=a^{m^iw'}$ (it follows from the fact that $\pi_1=<m^{\pi_1}>$), and it acts on an arbitrary element $a^v=a^{m^jv'}\in \widehat{Q}$, where $v'\in\pi_1'$.
\begin{align*}
\Phi((a^v)^x)&=\Phi(a^{m^jv'})*\Phi(q)=(m^{v'},v')*(m^{w'},w')\\
&=(m^{v'm^{w'}}, m^{-1}v'(w')^{-1}mw')=(m^{v'm^{w'}}, m^{-\epsilon(x)}v'(w')^{-1}mw')\\
&=(m^{v'},v')^{m^{w'}}=(m^{v'},v')^x=(\Phi(a^v))^x.
\end{align*}
We notice that $\Phi\Psi(m,1)=(m,1)$ and $\Psi\Phi(a)=a$. From the equivariance of the maps $\Phi$ and $\Psi$, and the fact that $\pi_1$ acts transitively on both quandles, follows that $\Phi\Psi=Id_{\tilde{Q}(\pi_1,m)}$ and $\Psi\Phi=Id_{\widehat{Q}}$.
Therefore, $\Phi\colon\widehat{Q}\to \tilde{Q}(\pi_1,m)$ is an isomorphism.
\end{proof}

\begin{figure}
\begin{center}
\includegraphics[height=4.5cm]{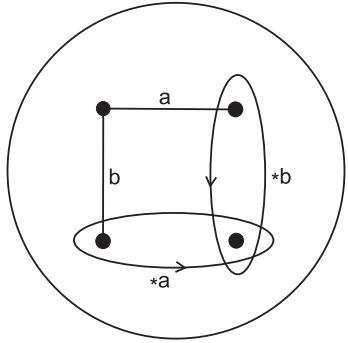}
\caption{Generators of ``the fundamental quandle of the closed trefoil knot."\label{closedquandle} }
\end{center}
\end{figure}

From the above analysis it is clear that if we allow $\delta_0$ to be rotated, we obtain the quandle of the closed trefoil knot.

\begin{corollary}
The fundamental quandle of the closed trefoil knot is isomorphic to the cord quandle on a 2-sphere with four punctures, generated by two arcs, $a$ and $b$,
as shown in the Figure \ref{closedquandle}. 
\end{corollary}

\end{document}